\documentclass[10pt]{article}

plus3pt minus2pt \abovedisplayshortskip=12pt plus3pt minus2pt
\usepackage{graphics, amsfonts, amssymb, latexsym, amsmath, amsxtra}
\usepackage[dvips]{graphicx}

\def\bd{\begin{displaymath}}
\def\ed{\end{displaymath}}
\def\be{\begin{equation}}
\def\ee{\end{equation}}
\newtheorem{theorem}{Theorem}
\newtheorem{cor}{Corollary}

\newtheorem{prop}{Proposition}

\newcommand {\mc} {\mathcal}
\newcommand {\p} {\partial}

\begin{document}

\begin{center}
\textbf{\Large{The modular hierarchy of  the Toda lattice}}
\end{center}

\large \vskip 1cm

\centerline{Maria A. Agrotis , Pantelis A. Damianou}

\bigskip

\centerline{ Department of Mathematics and Statistics, University
of Cyprus} \centerline { P. O. Box 20537, 1678 Nicosia, Cyprus}

Email addresses: agrotis@ucy.ac.cy (corresponding author),
damianou@ucy.ac.cy

\begin{abstract}
The modular vector field plays an important role in the theory of
Poisson manifolds and is intimately connected with the Poisson
cohomology of the space. In this paper we investigate its
significance in the theory of integrable systems. We illustrate in
detail the case of the Toda lattice both in Flaschka  and natural
coordinates.

\end{abstract}


\medskip

\noindent Classification number: 53D17,  37J35, 70H06.

\section{Introduction} \label{intro}

The infinitesimal generator of the modular automorphism group is
important in the cohomology theory of Poisson manifolds. On the
other hand, one of the main characteristics of integrable
Hamiltonian systems is the formation of hierarchies consisting of
Poisson tensors. Each one of these tensors carries naturally a
modular vector field. One may, therefore, speak of a hierarchy of
modular vector fields. The purpose of this paper is to begin an
investigation of the role of this sequence of vector fields in the
theory of integrable systems. As an illustrative example, we
choose the famous Toda lattice, a system that is one of the most
basic in the theory of finite dimensional integrable models and
closely related to the theory of simple Lie groups. There is no
doubt, however, that the results of this paper hold for other
similar systems as well.

The modular vector field appears first in a paper of Koszul
\cite{koszul}, as a special case of an operator on contravariant
differential forms of degree $-1$. The same vector field was later
used by Dufur and Haraki in \cite{dufour} to classify quadratic
Poisson brackets in $\mathbb{R}^3$. It was called "curl", which in
$\mathbb{R}^3$ is dual to "divergence".   In fact, the operator of
Koszul applied to a vector field, gives the usual divergence of
differential calculus with respect to the standard volume form.
One can therefore speak of the divergence of a Poisson bracket;
this is what the modular vector field is. The divergence of a
Poisson bracket is useful in the classification of Poisson
structures in low dimensions, e.g., \cite{xu,marmo}.

One class of operators considered by Koszul is the following. Let
$A_p$ be the space of covariant antisymmetric tensors on a smooth
 orientable  manifold $M$ of dimension $n$ and let $A^p$ denote the space of
covariant differential forms. Choose a volume form $\omega$ and
define an operator \bd D: A_p \to A_{p-1} \ed which satisfies \bd
D^2 =0 . \ed The volume form $\omega$ induces an isomorphism
$\Phi$ from $A_p$ to $A^{n-p}$.  Define the modular operator as:
\bd D=\Phi^{-1} \circ d \circ \Phi , \ed where $d$ is the exterior
derivative.

Weinstein \cite{weinstein2} gives a different interpretation.
Suppose one has a Poisson tensor $\pi$ and a smooth positive
volume form $\omega$. Consider the operator $X_{\omega},$ which
sends a smooth function $f$ to  ${\rm div}_{\omega} \mc{X}_f$,
where $\mc{X}_f$ is the Hamiltonian vector field generated by $f$
with respect to $\pi$. It turns out that this operator is a
derivation and hence a vector field. It coincides with the Koszul
operator, $D,$ at the level of contravariant 2--tensors. It
satisfies:
\begin{equation}
L_{X_{\omega} } \omega =0,   \quad  L_{X_{\omega}}
\pi=0.\label{aa1}
\end{equation}

Even though we deal with spaces that are orientable, we remark
that in the non--orientable case one can either replace volumes by
densities or use the recent approach of \cite{grab}.

 In local
coordinates $(x_1, \dots, x_n),$ the modular vector field of the
Poisson tensor $\pi$ is given by the formula:

\begin{equation} D(\pi)=\sum_{j=1}^{n}\left( \sum_{i=1}^n \frac{\partial
\pi_{ji}}{\partial x_i} \right)\frac{\partial}{\partial x_j}.
\label{divergence}\end{equation}

Lichnerowicz \cite{lich} considers the following cohomology
defined on the space of contravariant tensors of a Poisson
manifold. Let $(M, \pi)$ be a Poisson manifold.  Define a
coboundary operator $\partial_{ \pi}$ which assigns to each
$p$-tensor $A$, a $(p+1)$-tensor $\partial_{\pi} A$ given by \bd
\partial_{\pi}A=-[\pi, A] ,
\ed  where $[\cdot,\cdot]$ denotes the Schouten bracket. We have
that $\partial_{\pi}^2 A=[\pi, [\pi,A]]=0$ and consequently
$\partial_{\pi}$ defines a cohomology. More details on the
Schouten bracket can be found in \cite{lich,marsden,vaisman}. An
element $A$ is a $p$-cocycle if $[\pi, A]=0$. An element $B$ is a
$p$-coboundary if $B=[\pi,C]$, for some $(p-1)$-tensor $C$. Let
\bd Z^n(M, \pi) = \{ A \, : \,  [\pi, A]=0  \} \ed and \bd B^n(M,
\pi)= \{ B \,  : \, B=[\pi, C] \} \ . \ed The quotient \bd H^n(M,
\pi)= \frac{  Z^n(M, \pi)} {B^n(M, \pi) } \ed is the $n$th
cohomology group.  The elements of the first cohomology group are
the infinitesimal automorphisms of the Poisson structure modulo
the Hamiltonian vector fields. It follows from (\ref{aa1}) that
the modular vector field is an element of the first cohomology
group. If one replaces $\omega$ by $a \omega,$ where $a$ is a
positive smooth function, then

\be X_{a\omega} = X_{\omega} + \mc{X}_{\ln a} \label{aa3}  \ . \ee

\noindent Therefore, the modular vector field is a well--defined
element of the first cohomology group in the Lichnerowicz
cohomology. Since it forms an element of the cohomology group,
Weinstein  in \cite{weinstein2}  uses the term ''modular class''.
It makes sense to consider two such vector fields as equal if in
fact they differ by a Hamiltonian vector field. The reader can
also refer to \cite{lich,weinstein,vaisman} for more details on
Poisson manifolds and cohomology.

The modular operator $D$ is a graded derivation. The following
relations hold for a general 2-tensor $\pi$ and vector fields
$X,Y$:

\begin{eqnarray}
&& D[\pi,X]=[D(\pi),X]-[\pi,D(X)] \label{l1}\\
&& D[X,Y]=X D(Y)-Y D(X) .\label{l2}\end{eqnarray}  We note the
minus sign in (\ref{l1}) that appears because we use the
convention $\mc{X}_{f}^{\pi}(g)=\{f,g\}_{\pi},$ in the definition
of the Hamiltonian vector field $\mc{X}_{f}^{\pi}.$

The purpose of this paper is to present a complete study of the
modular hierarchy in the case of the Toda lattice. In Section
\ref{toda} we review the finite, classical Toda lattice and its
multi-Hamiltonian nature. In Section \ref{ab} we examine the
modular hierarchy  in Flaschka coordinates $(a,b) \in
\mathbb{R}^{2N-1}.$ In particular, we establish the Hamiltonian
character of the modular vector fields associated to the
well-known hierarchy of Poisson tensors for the Toda lattice, see
Theorem \ref{main2}, and present some of their basic properties in
Corollary \ref{cor1}. A new bi-Hamiltonian formulation of the Toda
lattice is provided in Theorem \ref{ham}. In Section \ref{qp} we
study the infinite modular sequence in natural coordinates $(q,p)
\in \mathbb{R}^{2N},$ and present a formula that iteratively
produces all members of this sequence. In Section \ref{summary} we
comment on the results of this paper and present them in compact
form.

\section{ The Toda lattice} \label{toda}

The Hamiltonian of the Toda lattice is given by  \be H(q_1, \dots,
q_N, \,  p_1, \dots, p_N) = \sum_{i=1}^N \, \frac{ 1 }{ 2} \,
p_i^2 + \sum _{i=1}^{N-1} \,  e^{ q_i-q_{i+1}} \ . \label{a1} \ee

\noindent This type of Hamiltonian was first considered by
Morikazu Toda \cite{toda}.  Equation (\ref{a1}) is known as the
classical, finite, non--periodic Toda lattice to distinguish the
system from the many and  various other versions, e.g.,  the
relativistic, quantum, infinite,  periodic etc.  The integrability
of the system was established in 1974 independently by Flaschka
\cite{flaschka1}, H\'enon  \cite{henon} and Manakov
\cite{manakov}.
 The original Toda
lattice can be viewed as a discrete version of the Korteweg--de
Vries equation. It is called a lattice, as in atomic lattice,
since interatomic interaction was studied. This  system also
appears  in cosmology and the work of Seiberg and Witten on
supersymmetric Yang--Mills theories. It has applications in analog
computing and numerical computation of eigenvalues. However, the
Toda lattice is mainly a theoretical mathematical model which is
important due to the rich mathematical structure encoded in it.

Hamilton's equations take the form

\bd
\begin{array}{lcl}
\dot q_j =p_j    \\
\dot p_j=e^{ q_{j-1}-q_j }- e^{q_j- q_{j+1}}, \quad j=1,\ldots, N.
\end{array}
\ed

\smallskip
\noindent The system is  integrable. One can find a set of
independent functions $\{  H_1, \dots,  H_N \} $  which are
constants of motion for Hamilton's equations. To  determine the
constants of motion, one uses Flaschka's transformation:
\begin{eqnarray}
&&  a_i  = \frac{1 }{ 2} e^{ \frac{1 }{ 2} (q_i - q_{i+1} ) },
\quad i=1,\ldots,N-1 \nonumber \\
&&  b_i  = -\frac{ 1 }{ 2} p_i, \hskip 1.5cm i=1,\ldots,N.
\label{a2} \
\end{eqnarray}

\smallskip
\noindent The equations of motion become \be
\begin{array}{lcl}
 \dot a _i& = & a_i \,  (b_{i+1} -b_i )    \\
   \dot b _i &= & 2 \, ( a_i^2 - a_{i-1}^2 ) .  \label{a3}
\end{array}
\ee

\smallskip
\noindent These equations can be written as a Lax pair  $\dot L =
[B, L] $, where $L$ is the Jacobi matrix

\bd \nonumber
 L= \begin{pmatrix} b_1 &  a_1 & 0 & \cdots & \cdots & 0 \cr
                   a_1 & b_2 & a_2 & \cdots &    & \vdots \cr
                   0 & a_2 & b_3 & \ddots &  &  \cr
                   \vdots & & \ddots & \ddots & & \vdots \cr
                   \vdots & & & \ddots & \ddots & a_{N-1} \cr
                   0 & \cdots & & \cdots & a_{N-1} & b_N   \cr
                   \end{pmatrix},
\ed

\smallskip
\noindent and

\bd \nonumber
    B =  \begin{pmatrix} 0 & a_1 & 0 & \cdots & \cdots &  0 \cr
                 -a_1 & 0 & a_2 & \cdots & & \vdots  \cr
                    0  & -a_2 & 0 & \ddots &  & \cr
                    \vdots &  & \ddots & \ddots & \ddots & \vdots \cr
                     \vdots & & &  \ddots & \ddots & a_{N-1} \cr
                     0 & \cdots &\cdots &  & -a_{N-1}  & 0 \cr
                     \end{pmatrix}.
\ed

\medskip \noindent This is an example of an isospectral
deformation; the entries of $L$ vary over time but the eigenvalues
remain constant, i.e. $\dot \lambda_i=0$. It follows that the
functions $ H_j=\frac{1 }{ j} {\rm tr} \, L^j$ are  constants of
motion. We note that \bd H_1=\lambda_1+ \lambda_2+ \dots
+\lambda_N=\sum_{i=1}^N b_i  \ed corresponds to the total momentum
and \bd H_2=\frac{1 }{ 2} \left( \lambda_1^2+\lambda_2^2+ \dots +
\lambda_N^2 \right) = \frac{1}{ 2} \sum_{i=1}^N b_i^2+
\sum_{i=1}^{N-1} a_i^2 \ed is the Hamiltonian.

Consider $\mathbb{R}^{2N} $  with  coordinates $(q_1, \dots , q_N,
p_1, \dots, p_N)$, the standard symplectic bracket \be \{ f, g
\}_s = \sum_{i=1}^N \left( \frac{\partial f }{ \partial q_i}
\frac{\partial g }{ \partial p_i} - \frac{\partial f }{ \partial
p_i} \frac{
\partial g }{
\partial q_i} \right) \nonumber,
\ee

\noindent and the mapping $F: \mathbb{R}^{2N} \to
\mathbb{R}^{2N-1}$ defined by \bd
 F:  (q_1, \dots, q_N, p_1, \dots, p_N) \to (a_1,  \dots, a_{N-1}, b_1, \dots, b_N) .
\ed

The Flaschka transformation $F$ is a symplectic realization of a
degenerate  Lie Poisson bracket  on  $\mathbb{R}^{2N-1}$, i.e.
there exists a Poisson bracket on $\mathbb{R}^{2N-1}$ which
satisfies \bd \{f, g \} \circ F = \{ f \circ F, g \circ F \}_s .
\ed

\noindent This  bracket  (up to a constant multiple) is given by
\be
\begin{array}{lcl}
\{a_i, b_i \}& =&-a_i  \\
\{a_i, b_{i+1} \} &=& a_i  \label{a4}   \ ;
\end{array}
\ee all other brackets are zero. We denote this bracket by
$\pi_1.$  Its Lie algebraic interpretation can be found in
\cite{kostant}.  In $\pi_1,$ the only casimir is $H_1=b_1+b_2 +
\dots +b_N,$ and the Hamiltonian is $H_2 = \frac{ 1 }{ 2}\ { \rm
tr}\ L^2$. The invariants $H_i$ are in involution with respect to
$\pi_1$. For a proof of these facts see \cite{damianou2}.

The quadratic Toda bracket appears in conjunction with isospectral
deformations of Jacobi matrices. Let $\lambda $ be an eigenvalue
of $L$ with normalized eigenvector $v$. Using standard
perturbation theory one obtains,

\bd \nabla \lambda = (2 v_1 v_2, \dots, 2  v_{N-1} v_N, v_1^2,
\dots, v_N^2)^t :=U^{\lambda}  , \ed

\noindent where $\nabla  \lambda $ denotes $( \frac{\partial
\lambda}{\partial a_1}, \ldots,\frac{\partial \lambda}{\partial
a_{N-1}},\frac{\partial \lambda
 }{ \partial b_1},\ldots,\frac{\partial \lambda
 }{ \partial b_N})$. Some manipulations show that $U^{\lambda}$ satisfies

\bd \pi_2 \,  U^{\lambda} = \lambda \,  \pi_1 \, U^{\lambda} , \ed
where $\pi_1$ and $\pi_2$ are skew-symmetric matrices. The
defining relations for the Poisson tensor $\pi_2$ are the
following quadratic functions of the $a_i$ and $b_i:$

\be \nonumber
\begin{array}{lcl}
\{a_i, a_{i+1} \}&=&\frac{ 1 }{ 2} a_i a_{i+1} \\
\{a_i, b_i \} &=& -a_i b_i                    \\
\{a_i, b_{i+1} \}&=& a_i b_{i+1}    \\
\{b_i, b_{i+1} \}&=& 2\, a_i^2   \ ;
\end{array}
\ee all other brackets are zero. The quadratic Toda bracket
appeared in a paper of Adler \cite{adler} in 1979. It is a Poisson
bracket in which the Hamiltonian vector field generated by $H_1$
is the same as the Hamiltonian vector field generated by $H_2$
with respect to the $\pi_1$ bracket. The tensor $\pi_2$ has ${\rm
det}\, L$ as casimir and $H_1 ={\rm tr}\, L$ as the Hamiltonian.
The eigenvalues of $L$ (and therefore  the $H_i$ as well) are
still in involution. Furthermore, $\pi_2$ is compatible with
$\pi_1$. We have \be \pi_2 \nabla H_j = \pi_1 \nabla H_{j+1}, \ \
j=1,2, \dots   \ . \label{a5} \ee These relations are similar to
the Lenard relations for  the KdV equation; they are generally
called the Lenard relations. Taking $j=1$ in (\ref{a5}), we
conclude  that the Toda lattice is bi--Hamiltonian.
Bi--Hamiltonian structures were introduced by Magri in
\cite{magri}. Using results from \cite{damianou5}, ones proves
that the Toda lattice is multi--Hamiltonian: \be \pi_2 \, \nabla
H_1=\pi_1 \, \nabla H_2 =\pi_0 \,  \nabla H_3 =\pi_{-1} \, \nabla
H_4= \dots \:. \label{a6} \ee

The Hamiltonian hierarchies  of the Toda lattice are  well-known.
The results are usually presented either in the natural $(q,p)$
coordinates or in the more convenient Flaschka coordinates
$(a,b)$. In the former case the hierarchy of higher invariants are
generated by the use of a recursion operator \cite{das,fernandes}.
We remark that recursion operators were first introduced by Olver
\cite{olver1}. The system is bi-Hamiltonian and one of the
brackets is symplectic. Thus, one can find a recursion operator by
inverting the symplectic  tensor. The recursion operator is then
applied to the initial symplectic bracket  to produce an infinite
sequence of Poisson tensors. However, in the case of the Toda
lattice in Flaschka variables $(a,b),$ the first two Poisson
brackets $\pi_1$ and $\pi_2$ are non-invertible and therefore this
method fails. The absence of a recursion operator for the finite
Toda lattice is also mentioned in Morosi and Tondo \cite{morosi},
where a Ninjenhuis tensor for the infinite Toda lattice is
calculated. The family of Poisson tensors in this case is
constructed using master symmetries. Invariant functions and
Hamiltonian vector fields are preserved by master symmetries. New
Poisson brackets are generated using Lie derivatives in the
direction of these vector fields, and they satisfy interesting
deformation relations. We quote the results from references
\cite{damianou1,damianou2}.

\begin{theorem} \label{lenard}
There exists  a sequence of vector fields $X_i$, for $i \ge -1$,
and a sequence of contravariant 2-tensors $\pi_j$, $j\ge 1$,
satisfying :

\smallskip
\noindent {\it i) } $\pi_j$   are all Poisson.

\smallskip
\noindent {\it ii) } The functions $H_i$,        $i\ge 1$, are in
involution
 with respect to all of the $\pi_j$.

 \smallskip
 \noindent
 {\it iii)}  $X_i (H_j) =(i+j) H_{i+j} $ ,  $i\ge -1$, $j\ge 1$.

 \smallskip
 \noindent
{\it iv)} $L_{X_i} \pi_j =(j-i-2) \pi_{i+j} $,  $i\ge -1$, $j\ge
1$.

\smallskip
\noindent {\it v)} $ [X_i, \ X_j]=(j-i)X_{i+j}$, $i\ge 0$, $j\ge
0$.

\smallskip
\noindent {\it vi)} $\pi_j\  \nabla \  H_i =\pi_{j-1}\   \nabla \
H_{i+1} $, where $\pi_j$ denotes  the Poisson matrix  of the
tensor $\pi_j$.
\end{theorem}

\noindent {Remark 1:} Theorem 3 was extended for all integer
values of the index in \cite{damianou5}.

\section{The modular class in Flaschka coordinates} \label{ab}

\noindent We consider the modular vector field of the Poisson
tensor $\pi_j,$ denoted by $Y_j=D(\pi_j),$ where $D$ is the Koszul
operator. We begin with some preliminary results needed to prove
the main theorem that establishes the Hamiltonian character of the
modular class. The vector fields of Theorem \ref{lenard}, denoted
by $X_j,$ are master symmetries for the Toda lattice system. For
example, $X_1$ is given as follows:

\begin{equation*}
X_1=\sum_{i=1}^{N-1} A_i \frac{\partial}{\partial a_i}+
\sum_{i=1}^{N} B_i \frac{\partial}{\partial b_i},
\end{equation*} where
\begin{eqnarray*}
&& A_i=-i a_i b_i + (i+2) a_i b_{i+1}    \\
&& B_i=(2i+3) a_i^2+(1-2i)a_{i-1}^2+b_i^2.
\end{eqnarray*}
If we define the function $f=\ln(a_1\cdots a_{N-1}),$ then we have
the following proposition:

\begin{prop} \label{prop1}

\begin{equation*}
D(X_1)=X_1(f)+2 H_1.
\end{equation*}

\end{prop}

\noindent \emph{Proof.}
\medskip

\noindent By definition
\begin{eqnarray*}
D(X_1)&=&\sum_{i=1}^{N-1}  \frac{\partial {A}_i}{\partial
a_i}+ \sum_{i=1}^{N}  \frac{\partial {B}_i}{\partial b_i}\\
&=& \sum_{i=1}^{n-1} [-ib_i+(i+2)b_{i+1}] + 2 \sum_{i=1}^{n} b_i
\\
&=& X_1(f)+2 \: \mbox{Tr}(L)=X_1(f)+2 H_1.
\end{eqnarray*}

\noindent A similar relation holds for the second master symmetry
$X_2$ that is defined as

\begin{equation*}
X_2=\sum_{i=1}^{N-1} C_i \frac{\partial}{\partial a_i}+
\sum_{i=1}^{N} D_i \frac{\partial}{\partial b_i},
\end{equation*} where
\begin{eqnarray*}
C_i &=&(2-i) a_{i-1}^2 a_i + (1-i) a_i b_i^2 + a_i b_i b_{i+1} \\
&+&(i+1) a_i a_{i+1}^2 +(i+1) a_i b_{i+1}^2  + a_i^3
+ (\sum_{j=1}^{i-1} b_j) a_i(b_{i+1}-b_i) \\
D_i &=& 2 (\sum_{j=1}^{i-1} b_j) a_i^2- 2 (\sum_{j=1}^{i-2} b_j)
a_{i-1}^2 + (2i+2) a_i^2 b_i + (2i+1) a_i^2 b_{i+1} \\ &+& (3-2i)
a_{i-1}^2 b_{i-1} +(4-2i) a_{i-1}^2  + b_i^3 .
\end{eqnarray*}

\begin{prop} \label{prop2}

\begin{equation*}
D(X_2)=X_2(f)+6 H_2.
\end{equation*}

\end{prop}

\noindent \emph{Proof.}

\medskip

\noindent Similar to the proof of Proposition \ref{prop1}.

\medskip

\noindent We can generalize the results of Propositions
\ref{prop1} and \ref{prop2} as follows:

\begin{prop}

\begin{equation*}
D(X_n)=X_n(f)+n(n+1) H_n.
\end{equation*}
\end{prop}

\noindent \emph{Proof.}

\medskip

\noindent The proof is inductive.  The result holds for $n=1,2$.
We assume that the result holds for $n \geq 2$ and prove it for
$n+1.$

\begin{eqnarray*}
D(X_{n+1})&=& \frac{1}{n-1} D[X_1,X_n]\\
&=& \frac{1}{n-1} \left( X_1D(X_n)-X_n D(X_1)\right)\\
&=& \frac{1}{n-1} \left( X_1(X_n(f)+n(n+1)H_n)-X_n (X_1(f)+2H_1)\right)\\
&=& \frac{1}{n-1} \left( [X_1,X_n]+n(n+1)X_1(H_n)-2X_1(H_n)\right) \\
&=& \frac{1}{n-1} \left( (n-1)X_{n+1}(f)+(n(n+1)-2)X_1(H_n)\right)\\
&=& X_{n+1}(f)+ (n+2)X_1(H_n)=X_{n+1}(f)+ (n+1)(n+2) H_{n+1}.
\end{eqnarray*}

Below we prove a statement that we will use in Theorem
\ref{main2}, but it is also important in its own right. It is a
new bi--Hamiltonian formulation of the Toda lattice.

\begin{theorem} \label{ham}
\begin{equation*}
\mc{X}_{H_1}^{\pi_j}=\mc{X}_{g}^{\pi_{j+1}}, \:\:\: where\:\:\:
g=\ln(\det(L)), \:\: and \:\: j\geq 1.\end{equation*}
\end{theorem}

\noindent \emph{Proof.}

\medskip

\noindent We use the Lenard relations for the eigenvalues
\cite{damianou2,damianou6}, $\pi_j \nabla \, \lambda_i= \lambda_i
\: \pi_{j-1} \nabla \, \lambda_i .$
\begin{eqnarray*}
\mc{X}_{H_1}^{\pi_j}&=&\mc{X}_{(\lambda_1 + \dots +\lambda_N)}^{\pi_j} \\
             &=& \pi_j \nabla  (\lambda_1 + \dots +\lambda_N) \\
             &=& \sum_{i=1}^N \pi_j \nabla \lambda_i  .
\end{eqnarray*}

\noindent On the other hand
\begin{eqnarray*}
\mc{X}_g^{\pi_{j+1}}&=&\mc{X}_{\ln{\lambda_1 \dots \lambda_N}}^{\pi_{j+1}} \\
             &=& \pi_{j+1} \nabla \ln{\lambda_1 \dots \lambda_N }\\
             &=& \sum_{i=1}^N \pi_{j+1}\frac{1 }{ \lambda_i} \nabla
             \lambda_i \\
             &=& \sum_{i=1}^N \pi_j \nabla \lambda_i \\
             &=& \mc{X}_{H_1}^{\pi_j} .
\end{eqnarray*}

The following theorem investigates the Hamiltonian character of
the divergence of the infinite sequence of Poisson tensors
$\pi_j.$ In particular, it states that the divergence of the
Poisson tensor $\pi_j$ is the Hamiltonian vector field given by
the function $h:=\ln(a_1 \cdots a_{N-1})+ (j-1) \ln(\det(L)).$

\begin{theorem} \label{main2}

$$Y_j=\mc{X}_{f+(j-1) \: g}^{\pi_j},\:\:
where \:\: f=\ln(a_1 \cdots a_{N-1}),\:\: and
\:\:g=\ln(\det(L)).$$
\end{theorem}

\noindent \emph{Proof.}

\medskip

\noindent We will prove the proposition inductively in two steps.
First we will show that it holds for $j=1$ and then for $2 \leq j
\leq 3.$ Consequently, we will show that it holds for $j=4$ and
then for $j \geq 5.$

We have that $Y_1=D(\pi_1),$ where $\pi_1$ is given by
\begin{eqnarray*}
&& \{a_i,b_i\}=-a_i \\
&& \{a_i,b_{i+1}\}=-a_i.
\end{eqnarray*}

\noindent Thus, using (\ref{divergence}) we obtain
$$Y_j=\frac{\partial}{\partial b_1} - \frac{\p}{\p
b_n}=(0,\ldots,0,1,0,\ldots,0,-1)^t.$$ On the other hand, it is
not hard to show that the Hamiltonian vector field induced by
$f=\ln(a_1\cdots a_{N-1})$ with respect to the bracket $\pi_1$ has
the form:
\begin{eqnarray*}
&& \dot{a}_i=0, \quad i=1, \ldots, N-1 \\
&& \dot{b}_1=1, \quad \dot{b}_N=-1 \\
&& \dot{b}_i=0, \quad i=2,\ldots, N-1.
\end{eqnarray*} Hence $Y_1=\mc{X}_{f}^{\pi_1}.$

Now, let us assume that the theorem holds for $j \geq 4 .$ We will
then show that it holds for $j+1.$ By definition
$L_{X_1}\pi_j=[X_1,\pi_j]=[\pi_j,X_1]$ and using $(iv)$ of Theorem
\ref{lenard}, we have that $\pi_{j+1}=\frac{1}{(j-3)} [\pi_j,
X_1].$ Thus,

\begin{eqnarray*}
Y_{j+1} &=& D(\pi_{j+1}) = \frac{1}{(j-3)} D([\pi_j,X_1]) \\
& = & \frac{1}{(j-3)} \left( [D(\pi_j),X_1]- [\pi_j,D(X_1)]\right) \\
& = & \frac{1}{(j-3)} \left( [Y_j,X_1]- [\pi_j,D(X_1)]  \right) \\
& = & \frac{1}{(j-3)} \left( [\mc{X}_{f}^{\pi_j},X_1]+(j-1)[\mc{X}_{g}^{\pi_{j}},X_1]- [\pi_j,D(X_1)]  \right) \\
& = & \frac{1}{(j-3)} \left( [[f,\pi_j],X_1]+(j-1)[\mc{X}_{H_1}^{\pi_{j-1}},X_1]- [\pi_j,D(X_1)]  \right) \\
& = & \frac{1}{(j-3)} \left( [\pi_j,[X_1,f]]+[f,[\pi_j,X_1]]+(j-1)(j-2)\mc{X}_{H_1}^{\pi_{j}}- [\pi_j,D(X_1)]  \right) \\
& = & \frac{1}{(j-3)} \left( [\pi_j,X_1(f)]+(j-3) [f,\pi_{j+1}]+(j-1)(j-2)\mc{X}_{H_1}^{\pi_{j}}- [\pi_j,D(X_1)]  \right) \\
& = & \frac{1}{(j-3)} \left( [\pi_j,D(X_1)]-2 [\pi_j,H_1]+(j-3)
\mc{X}_{f}^{\pi_{j+1}}+(j-1)(j-2)\mc{X}_{H_1}^{\pi_{j}}- [\pi_j,D(X_1)]  \right) \\
& = & \mc{X}_{f}^{\pi_{j+1}}+\frac{(j-1)(j-2)-2}{(j-3)}
\mc{X}_{H_1}^{\pi_j}\\
& = & \mc{X}_{f}^{\pi_{j+1}}+j \mc{X}_{g}^{\pi_{j+1}}.
\end{eqnarray*}
In an identical manner one proves that the theorem holds for $2
\leq j \leq 3.$ To show that the proposition holds for $j=4$ we
use a similar argument as the one used above, however, we employ
the relation
$\pi_4=-\frac{1}{2}L_{X_2}\pi_2=-\frac{1}{2}[X_2,\pi_2],$ to avoid
division by zero.

\begin{cor}  \label{cor1}

\hskip 1cm \vskip .3cm

$a) \:  Y_i (H_j)=Y_j(H_i) $

$b) \: L_{Y_i} \pi_j=-L_{Y_j} \pi_i.$
\end{cor}

\noindent \emph{Proof.}

\medskip

\noindent For part a) we have:
\begin{eqnarray*}
 Y_i (H_j) &=& \left( X_f^{\pi_i} + (i-1) \mc{X}_{H_1}^{\pi_{i-1}} \right) (H_j) \\
 \       &=& \mc{X}_f^{\pi_i} (H_j)  \\
  \    &=& \{f, H_j \}_{\pi_i} \\
        &=& \nabla f \pi_i \nabla H_j \\
          &=& \nabla f \pi_j \nabla H_i \\
           &=& \{f, H_i \}_{\pi_j} \\
           &=& Y_j (H_i).
           \end{eqnarray*}

\noindent For part b), \be \nonumber \left[ [h, \pi_j], \pi_i
\right] + \left[ [\pi_j, \pi_i], h \right] + \left[ [ \pi_i, h],
\pi_j \right]=0. \ee

\noindent Since $[\pi_i, \pi_j]=0$ the result follows.

\section{The modular hierarchy in $(q,p)-$coordinates} \label{qp}

The bi--Hamiltonian structure of the Toda lattice appears in a
paper of  Das and Okubo in 1989 \cite{das}.  The generation of
master symmetries and Poisson tensors using Oevel's theorem and
the connection with the results in Flaschka coordinates is due to
Fernandes \cite{fernandes}. In principle, the  method is general
and  works for other finite dimensional systems as well. For
example, this  approach was used  by Nunes da Costa and Marle
\cite{joana}  in the case of the relativistic Toda lattice. The
procedure goes as follows. One defines a second Poisson bracket in
the space of canonical variables $(q_1, \dots, q_N, p_1,\dots,
p_N)$. This gives rise to a recursion operator, call it $\mc{R}$.
The presence of a conformal symmetry, as defined by Oevel, allows
one to use the recursion operator and generate an infinite
sequence of master symmetries. These, in turn, project to the
space of the new variables $(a,b)$ to produce a sequence of master
symmetries in the reduced space. We quote the relevant theorem of
Oevel that appears in \cite{oevel2}.

\begin{theorem} \label{th2}
Suppose that   $Z_0$ is a conformal symmetry for both Poisson
tensors $J_1$, $J_2$ and  function $h_1$, i.e.  for some scalars
$\lambda$, $\mu$, and $\nu$ we have \bd {\cal L}_{Z_0} J_1=
\lambda J_1, \quad{\cal L}_{Z_0} J_2 = \mu J_2, \quad {\cal
L}_{Z_0} h_1 = \nu h_1  . \ed Then the vector fields $Z_i = {\cal
R}^i Z_0$ are master symmetries and we have,
\begin{eqnarray*}
&& (a) \ {\cal L}_{Z_i} h_k = (\nu +(k-1+i) (\mu -\lambda))
h_{i+k}
\\
&& (b) \ {\cal L}_{Z_i}J_k = (\mu +(k-i-2) (\mu -\lambda))
J_{i+k} \\
&& (c)  [Z_i, Z_k]= (\mu - \lambda) (k-i) Z_{i+k}  .
\end{eqnarray*}
\end{theorem} We note that master symmetries were introduced in
\cite{fokas}.

We proceed with the computation of the modular vector fields that
are associated to the family of Poisson tensors given in
$(q,p)-$coordinates. The first Poisson tensor in the hierarchy is
the standard canonical symplectic tensor: \bd
\hat{J}_1=\begin{pmatrix} 0 & I_{N}
\\-I_{N} & 0
\end{pmatrix},
\ed where $I_N$ denotes the $N \times N$ identity matrix. The
second Poisson tensor has the form,

\begin{equation*}
\hat{J}_2=\begin{pmatrix} A_{N} & B_{N} \\-B_{N} & C_{N}
\end{pmatrix},
\end{equation*} where $A_N$ is the $N \times N$ skew-symmetric matrix defined by
$a_{ij}=1=-a_{ji},\:\mbox{for} \: i<j,$
$B_{N}=\mbox{diag}(-p_1,-p_2,\ldots, -p_N),$ and $C_{N}$ is the $N
\times N$ skew-symmetric matrix whose non-zero terms are given by
$c_{i,i+1}=-c_{i+1,1}=e^{q_i-q_{i+1}},\: \mbox{for} \:
i=1,2,\ldots, N-1.$ If we let $J_1=4 \hat{J}_{1}$ and $J_2=2
\hat{J}_{2}$ then $J_1$ and $J_2$ are mapped precisely onto the
brackets $\pi_1$ and $\pi_2$ under the Flaschka transformation.

It is easy to see that we have a bi-Hamiltonian pair. We define
\bd h_1=-2(p_1+p_2+\dots +p_N)  , \ed and $h_2$ to be the
Hamiltonian \bd h_2=\sum_{i=1}^N   \frac{ 1 }{ 2}  p_i^2 + \sum
_{i=1}^{N-1} \,  e^{ q_i-q_{i+1}}   . \ed Under Flaschka's
transformation  (\ref{a2}),  $h_1$ is mapped onto $4(b_1+b_2+
\dots+b_N)=4 \  {\rm tr} L= 4 H_1$ and $h_2$ is mapped onto $2 \
{\rm tr} L^2=4 H_2$. Using  the relationship \bd \pi_2 \nabla \
H_1=\pi_1 \nabla \  H_2  , \ed which follows from part $(iv)$ of
Theorem \ref{lenard}, we obtain, after  multiplication by 4, the
following pair: \bd J_1 \nabla \ h_2= J_2 \nabla \  h_1  . \ed

\noindent The Lenard relations for the eigenvalues translate into
\be \pi_j \nabla \, \lambda_i= \lambda_i \: \pi_{j-1} \nabla \:
\lambda_i  . \label{eigenvalues} \ee

\noindent We define the recursion operator as follows: \bd {\cal
R}=J_2 J_1^{-1} . \ed

\noindent This operator raises degrees and we therefore   call it
the positive Toda operator. In $(q,p)$ coordinates, we use the
symbol $\mc{X}_i$ as a shorthand for  $\mc{X}_{h_i}^{J_1},$ the
Hamiltonian vector field of $h_i$ with respect to the symplectic
bracket $J_1$. It is generated, as usual, by

\bd \mc{X}_i = {\cal R}^{i-1} \mc{X}_1 . \ed In a similar fashion
we obtain the  higher order Poisson tensors \bd J_i = {\cal
R}^{i-1} J_1 . \ed

\noindent We then define the conformal symmetry \bd
Z_0=\sum_{i=1}^N (N-2i+1) \frac{\partial }{ \partial q_i}
+\sum_{i=1}^N p_i \frac{\partial }{ \partial p_i} . \ed

\noindent It is straightforward to verify that

\begin{eqnarray*}
&& {\cal L}_{Z_0} J_1=- J_1, \\
&&{\cal L}_{Z_0} J_2=0.
\end{eqnarray*}

\noindent In addition,
\begin{eqnarray*}
&& Z_0(h_1)=h_1, \\
&& Z_0(h_2)=2h_2  .
\end{eqnarray*}

\noindent Consequently, $Z_0$ is a conformal symmetry for $J_1$,
$J_2$ and $h_1$. The constants appearing in Theorem \ref{th2} are
$\lambda=-1$, $\mu=0$ and $\nu=1$. Therefore, we end--up with the
following deformation relations:

\begin{equation*}
[Z_i, h_k]= (i+k)h_{i+k}
\end{equation*}

\begin{equation*}
L_{Z_i}  J_k = (k-i-2) J_{i+k}
\end{equation*}

\begin{equation*}
[ Z_i, Z_k ]  = (k-i) Z_{i+k}.
\end{equation*}

We compute the divergence of the master symmetry $Z_1,$ as we will
use it in the proof of Theorem \ref{main3}. It was proved in
\cite{damianou1} that $Z_1=\sum_{i=1}^{N} \lambda_i^2 \frac{\p}{\p
\lambda_i},$ where $\lambda_i$ are the eigenvalues of the Jacobi
matrix $L.$ Using the definition of divergence we have that
$D(Z_1)=2 \sum_{i=1}^{N}\lambda_i=2 \sum_{i=1}^{N}
b_i=-\sum_{i=1}^{N} p_i=\frac{1}{2}h_1.$ Therefore, \bd
D(Z_1)=\frac{1}{2}h_1.   \ed

The following proposition will also come to use in the proof of
Theorem \ref{main3}.

\begin{prop} \label{ham2}
\be \mc{X}_{h_1}^{J_j}=4 \mc{X}_f^{J_{j+1}}, \:\: where \:\:
f=\ln(\sqrt{\det \mc{R}})\:\: for \:\: j\ge 1 \: .\ee

\end{prop}

\noindent \emph{Proof.}

\medskip

\noindent Similar to the proof of Theorem \ref{ham}. We note that
one uses the fact that the eigenvalues of $\mc{R}$ are the squares
of the eigenvalues of $L$, see \cite{damianou6,falqui}.

\begin{theorem}  \label{main3} For $j \ge 1$,  $Y_j$ is a Hamiltonian
vector field given as  \begin{equation} Y_j=(j-1)
\mc{X}_{f}^{J_j}, \:\: where\:\: f= \ln(\sqrt{{\rm det}\, {\cal
R}})  .
\end{equation}
\end{theorem} We will prove the theorem for j=1, \ldots, 4, and then
we will use induction for $j \ge 5.$ First we observe that
$Y_1=\vec{0},$ since $J_1$ is symplectic. Using the general form
of the tensor $J_2$ we obtain the following: $Y_2=-2 \sum_{i=1}^N
\frac{ \p }{ \p q_i}.$ A simple calculation gives that
$\mc{X}_{h_1}^{J_1}=-8\sum_{i=1}^N \frac{ \p }{ \p q_i}.$ Thus
$Y_2=\frac{1}{4} \mc{X}_{h_1}^{J_1}=\mc{X}_{f}^{J_2},$ using
Proposition \ref{ham2}.


We have that $J_3=-[Z_1, J_2]$ and $D(Z_1)=\frac{1}{2}h_1.$
Therefore,

\begin{eqnarray*}
Y_3=D(J_3)& = & -D[J_2,Z_1] \\
&=& -\left([D(J_2), Z_1]-[J_2, D(Z_1)] \right) \\
&=& -[Y_2,Z_1]+\frac{1}{2} [J_2,h_1] \\
        &=& -[\mc{X}_f^{J_2}, Z_1]+\frac{1 }{ 2}
        \mc{X}_{h_1}^{J_2}.
\end{eqnarray*}

\noindent Using the super--Jacobi identity for the Schouten
bracket, the first term equals \bd -[J_2, [Z_1,f]-[f,[J_2,Z_1]].
\ed \noindent Therefore
\begin{eqnarray*}
Y_3 &=& -\frac{1}{4} \mc{X}_{h_1}^{J_2}+
[f,J_3] +\frac {1 }{ 2} \mc{X}_{h_1}^{J_2} \\
        &=& \frac{1}{4} \mc{X}_{h_1}^{J_2}+\mc{X}_f^{J_3}\\
        &=& 2 \mc{X}_f^{J_3}.
        \end{eqnarray*}

In the last step of the argument we have used that
$\mc{X}_{h_1}^{J_2}=4 \mc{X}_f^{J_3}$. The proof of the formula
$Y_{4}=3 \mc{X}_f^{J_4}$ is identical to the one for $Y_3$ except
that we employ the relation $[Z_2, J_2]=-2 J_4$. An inductive
argument based on the same technique that we have used for $Y_3,$
can also be used to show that the result of the theorem holds for
$j \ge 5.$ We omit the details.

In the theorem that follows, we present an iterative formula that
produces all members of the modular class in terms of the
recursion operator $\mc{R}$ and the modular vector field $Y_2=-2
\sum_{i=1}^N \frac{ \p }{ \p q_i}=(-2,\ldots,-2,0,\ldots,0)^t.$

\begin{theorem}
For $j \ge 2$, \bd Y_{j+1}=j \: \mc{R}^{j-1} \:Y_2 . \ed
\end{theorem}

\noindent \emph{Proof.}
\begin{eqnarray*}
Y_{j+1} &=& j \mc{X}_f^{J_{j+1}}  \\
        &=&\frac {j }{ 4} \mc{X}_{h_1}^{J_j}=\frac {j }{ 4} \mc{X}_{h_j}^{J_1}=\frac {j }{ 4} \mc{X}_{j} \\
        &=& \frac{j }{ 4} \mc{R}^{j-1} \mc{X}_{1} \\
        &=& \frac{j }{ 4} \mc{R}^{j-1} \mc{X}_{h_1}^{J_1} \\
        &=& \frac{j }{ 4} \mc{R}^{j-1} 4 \mc{X}_{f}^{J_2} \\
        &=& j \: \mc{R}^{j-1} \:Y_2.
\end{eqnarray*}

\section{Summary} \label{summary}

In this paper we study the hierarchy of modular vector fields
associated to the infinite family of Poisson tensors for the
classical Toda lattice equations. We present analytical
expressions of the modular vector fields both in Flaschka
coordinates $(a,b)\in \mathbb{R}^{2N-1}$, as well as in natural
coordinates $(q,p) \in \mathbb{R}^{2N}.$ In both cases, all the
members of the infinite modular family, denoted by $Y_j, \: j \geq
1,$ are Hamiltonian. In $(a,b)-$variables we have that

$$Y_j=\mc{X}_{\ln(a_1 \cdots a_{N-1})+(j-1)\ln(\det(L))}^{\pi_j}$$
where $L$ is the Jacobi matrix of the Lax pair. In natural
coordinates $(q,p),$ the modular vector field takes the form

$$Y_j=(j-1) \mc{X}_{\ln(\sqrt{\det{\mc{R}}})}^{J_j}$$ where
$\mc{R}$ is the recursion operator. It is not difficult to show
that the term $\mc{X}_{\ln(\det(L))}^{\pi_j}$ is the projection of
the vector field $\mc{X}_{\ln(\sqrt{\det{\mc{R}}})}^{J_j}$ under
the Flaschka map. The term $\mc{X}_{\ln(a_1 \cdots \:
a_{N-1})}^{\pi_j}$ makes its appearance due to the change of
coordinates.

The following properties are proved regarding the behavior of
modular vector fields when they are applied on the constants of
motion $H_j,$ and the Lie derivative of the Poisson tensor $\pi_j$
in the direction of the modular vector field $Y_i:$
\begin{eqnarray*}
&& i) \:  Y_i (H_j)=Y_j(H_i) \\
&& ii) \: L_{Y_i} \pi_j=-L_{Y_j} \pi_i.
\end{eqnarray*}

In $(q,p)-$variables, we presented a formula that iteratively
produces all members of the infinite family of modular vector
fields, in terms of the recursion operator. Namely,
\begin{eqnarray*}
&& Y_1=\vec{0}, \\
&& Y_2=-2 \sum_{i=1}^{N} \frac{\partial}{\partial
q_j}=(-2,\ldots,-2,0,\ldots,0)^t, \\
&& Y_{j+1}=j \: \mc{R}^{j-1} \:Y_2 \quad \mbox{for} \quad j \ge 2.
\end{eqnarray*}
We conclude with an alternate bi-Hamiltonian formulation of the
Toda lattice given by the relation
$\mc{X}_{H_1}^{\pi_j}=\mc{X}_{\ln(\det(L))}^{\pi_{j+1}},$ in
addition to the existing one
$\mc{X}_{H_1}^{\pi_{j+1}}=\mc{X}_{H_2}^{\pi_{j}}.$

\section{Acknowledgements}

The authors would like the thank the Cyprus Research Promotion
Foundation for support through the grant CRPF0504/03. We would
also like to thank G. Marmo for pointing out the importance of the
modular class in Poisson Geometry.

\end{document}